\documentclass [10pt]{article}
\usepackage[utf8]{inputenc}
\usepackage{amsmath}
\usepackage{amsfonts}
\usepackage{amssymb}
\usepackage{graphicx}
\usepackage{soul}
\usepackage{parskip}
\usepackage{natbib}
\usepackage{bm}
\usepackage{algorithm}
\usepackage[noend]{algpseudocode}
\usepackage{blindtext}
\usepackage{booktabs}
\usepackage{graphicx}
\usepackage{titlesec}

\title{\textbf{Structural Control Analysis of System Dynamics Models}} 
\author{
\textsc{Tianyi Li}\thanks{tianyil@mit.edu} \\[1ex] 
System Dynamics Group, Sloan School of Management, MIT \\ 
 \\
}
\date{} 

\begin{document}
\maketitle

\begin{abstract}
Structural control theory could be applied to study the control principles of social, economic and managerial systems. System Dynamics (SD) is the target field in social-economic sciences for endogenizing this theory, a subject that provides modeling solutions to real-world problems. SD models adopt diagrammatic representations, making it an ideal ground for transplanting structural control theory which utilizes similar graphic representations. This study sets up the theoretical ground for conducting structural control analysis (SCA) on SD models, summarized as a post-modeling workflow for SD practitioners, which serves as a specific application of the general structural control theory in social-economic sciences. Theoretical and practical establishments for SCA components are developed coordinately. Specifically, this study addresses the following questions: (1) How do SD models differ from physical control systems in graphic representations, and how do these differences affect the way of applying structural control theories to SD? (2) How could one identify control inputs in SD models, and how could different levels of system control in SD models be conceptualized? (3) What are the structural control properties for important SD components, and how could these properties and control principles help justify modeling heuristics in SD practice? (4) What are the procedures for conducting Structural Control Analysis (SCA) in SD models, and what are the implications of SCA results for model calibration and decision making? Overall, this study provides general insights for system control analysis of nonlinear dynamic simulation models, which may go beyond SD and extend to various disciplines in social-economic sciences.
\end{abstract}


Founded by Jay Forrester \citep[e.g.][]{F1970, F1997} and primarily promoted by John Sterman \citep[e.g.][]{S2000}, System Dynamics (SD) is the field that provides modeling solutions to real-world problems, typically social, managerial and organizational problems. SD views the underlying problem from a systematic point of view, highlighting interdependent relationships between model variables and \textit{feedback} structures of the system. In a typical SD model (Figure 1), one identifies \textit{stock variables} (variables in boxes) which represent level components of the system that are of major concerns; arrows indicate interdependencies between variables, including relationships that could be represented by normal functionals, and other nonlinear relationships that one may characterize using \textit{table functions} which graphically specify the desired dependencies. SD represents the type of modeling practice in social and managerial sciences that takes the system/network perspective, which methodologically differs from econometrical approaches such as regressions, game theories etc. It provides modelers and audiences with a straightforward system view of the problem to be solved, along with direct handles for policy decision-making and strategy analysis, although the quantitative rigor of its solution is slightly compromised compared with econometrical attempts. 

\begin{figure}[!ht]
\centering   
	\includegraphics[width=6in]{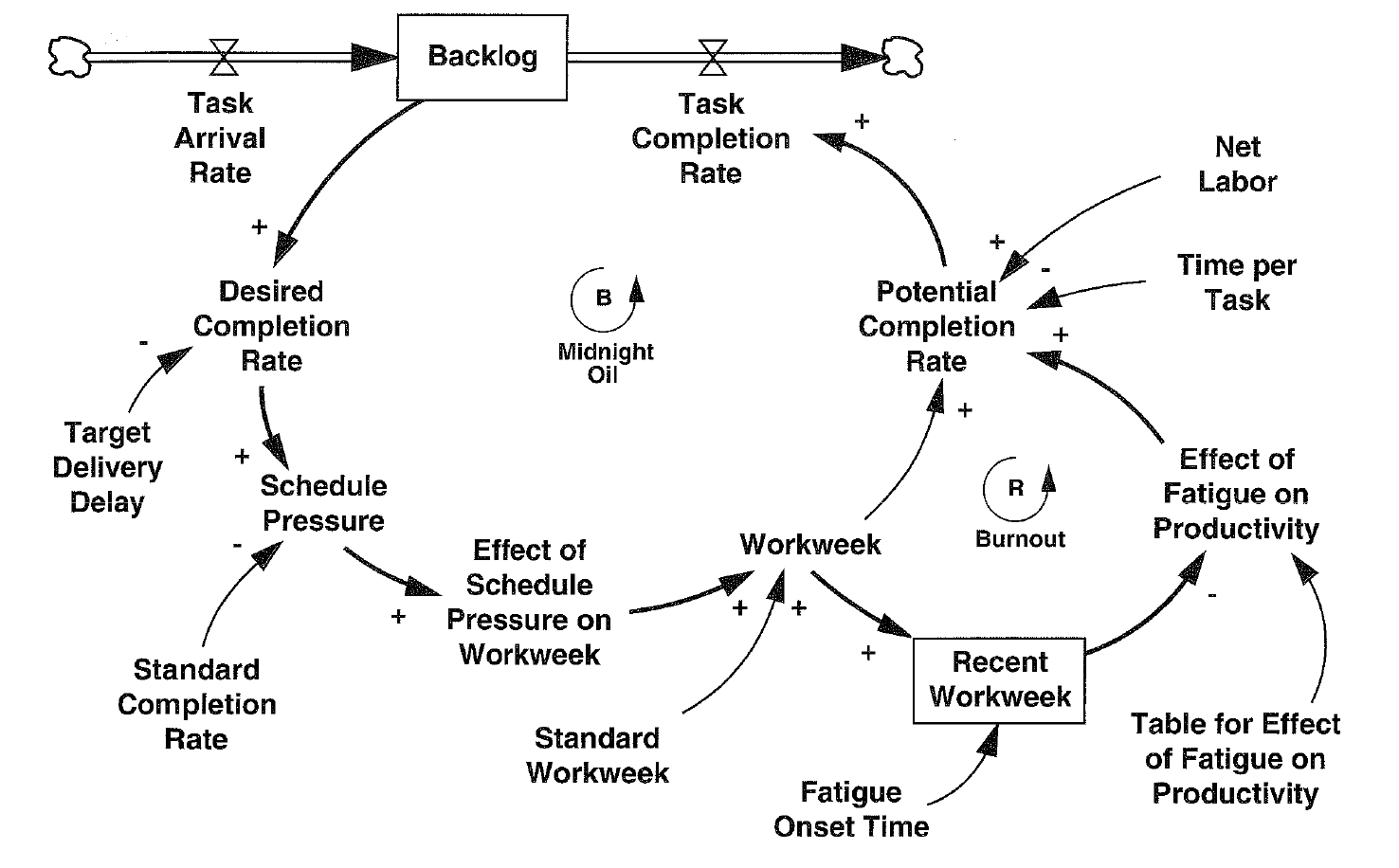}  
\caption{A typical SD model, showing stock and flow structures, table functions, balancing (B) and reinforcing (R) feedback loops. Adapted from Figure 14-13 in \citet{S2000}.} 

\end{figure}

Control theory constitutes the major theoretical ground for System Dynamics. As the basic theories in SD have been substantially developed by Forrester back in the 70s and 80s, the field has mainly been focusing on applications and pedagogical concerns in the past few decades. Nevertheless, over the years, the theoretical contents of SD have been steadily enlarged in the field, assimilating new advances in control theories as well as other related technical areas. Typical studies of this type include \citet{OS1989}, which demonstrated and studied the linearization process as well as the controller design of nonlinear economic systems, and \citet{Eet2014}, which adopted Hearn's perturbation method in the sensitivity analysis of table functions. Notably, a major stream of theoretical efforts has been spent on studying the structural dominance of SD models, and an abundant set of mathematical tools have been developed accordingly, including the pathway participation metrics (PPM) \citep{Met2004}, and more importantly, the eigenvalue analysis LEEA (loop eigenvalue elasticity analysis) and DWWA (dynamic decomposition weight analysis), which are successfully developed through a series of attempts \citep{Set2010,KO2006,KO2008,K2012,O2016,NO2018}. Inspired by those advances, it is expected that more and more theoretical inputs borrowing from other disciplines will take place in the SD field, help consolidate as well as amplify the theoretical ground of SD. With such goals in mind, this study serves as a tentative attempt.

\subsection*{\textit{Why do we need to study the control properties of SD models?}}

In SD models, as in typical dynamic systems as well, most variables are often not directly tunable and it is important that one could tune other variables (i.e., control inputs) to access them. Reference modes demonstrate the desired trajectories of variables, but one often does not know before completing the model whether he could drive the variables to reach certain values, i.e. whether the system is controllable. He also often has little idea whether the model is built with flaws in system control. A number of heuristics have been established for this concern, to check if the model has embodied inferior practice in its formulation. For example, it is generally acknowledged that stocks should ensure first-order exit control, that the same parameters should not take effect on multiple places, that table functions are a versatile tool for specifying nonlinear relationships, and that balancing feedbacks should be called for to stabilize the system. All of these heuristics point to the controllability of the system. On the opposite, when tuning the parameters after one finishes the model, certain phenomena are often signals of bad formulations: some stocks start to have negative values, or that some extreme values break down the simulation. All these undesired phenomena indicate that the model has drawbacks in system control and needs re-engineering. Given such concerns on both sides, it is then necessary to investigate the control principles for SD models and, if possible, secure certain theoretical grounds for some of those tested heuristics in practice. 

When discussing undesirable system behaviors in simulation, it is important to note that, to certain extent, the desired control properties of SD models are \textit{orthogonal to} that of physical dynamic systems. For physical systems, it is significant that the system is fully controllable, i.e., control inputs could always successfully manipulate the system according to human desires, and any uncontrolled part of the physical system poses a risk to human beings. However, for most SD models, the concern is the reverse: since SD models simulate dynamic processes in real life, in which model variables are often not of physical nature, it is strongly desired that the model could be operational, in the sense that the dynamics of real-life variables are confined to real-life trajectories, as opposed to incorporating pathways that are impossible in practice, i.e., in this case, the model should \textit{not} be fully controllable. Indeed, it is a common critique of SD models, if not all dynamic and simulation models applied in social, economic and managerial areas, that those models are too versatile, overly parameter-tuning-dominant, and essentially providing no insights to the underlying problem if ``anything can go''. It is expected that dynamic models in social sciences show that only a limited set of trajectories are possible for key variables, and even that models could generate natural pathways by themselves, exempt of exogenous inputs. As one could see, this concern highlights the importance of controllability analysis on SD models from a different angle. It is believed that, under established control principles, if an SD model (or a general simulation model in social sciences), could be shown to only hold \textit{limited controllability}, the model would then be able to general \textit{more} convincing power in its utilities in policy analysis and decision making.

Another significant concern with SD and general simulation models is the difficulty in designing calibration strategies for a nontrivial number of model variables. Large dynamic models may contain up to hundreds of variables, and it is an important issue to figure out which of these variables are of primary concern. To solve this problem, for a long time, people across disciplines have been designing techniques for the task of model dimensionality reduction, especially for non-linear systems \citep[e.g.,][]{Tet2000}. Generally speaking, variables that influence the system's dynamics most should have high priorities in data acquisition and model calibration. The characterization of a variable's influence on the model could be interpreted in different ways. For example, \citet{O2004} suggested multiple model partition methods for SD, which essentially highlighted variables' (and loops') structural importance on models. Inspired by this idea, the structural control analysis on SD models may offer an alternative guidance for the design of model calibration strategies. Given an SD model, modelers are able to inspect the control properties of each identified control input; different control inputs may have different effect on the model: some may exert full control, while in most case a single control input could only affect a certain part of the model. It is then possible to evaluate and rank different model control inputs by their significance on the system control, and by which means provide the modeler with guidances for model calibration strategies. Standardization of this evaluation process would be a useful post-modeling toolkit for SD practitioners.

\section*{Structural Controllability of Dynamic Systems}

The controllability of linear dynamic systems is an important topic in control theory \citep{LB2016}. A dynamic system is \textit{controllable}, if and only if one can drive the system from any initial state to any final state within finite time, by manipulating the control inputs. The best-known criterion for dynamic system's controllability is Kalman's rank condition \citep{K1963}. For a linear time-invariant (LTI) system $(\bm{A},\bm{B})$:
\begin{equation}
\dot{\bm{x}}(t) = \bm{A} \bm{x} (t) + \bm{B} \bm{u} (t)
\end{equation}

where $\bm{A} = N\times N$, $\bm{B} = N\times M$, $\bm{x}$ and $\bm{u}$ are state variables and control inputs, respectively. The system is controllable if and only if the controllability matrix:
\begin{equation}
\bm{C} \equiv [\bm{B}, \bm{AB}, \bm{A}^2\bm{B}, ...\bm{A}^{N-1}\bm{B}] 
\end{equation}
has full rank, i.e., $rank(\bm{C}) = N$.

Kalman's rank condition is rigorous and works on all kinds of dynamic systems. However, it is often impossible to apply this criterion in real world systems, since in many complex networks the system parameters (elements in $\bm{A}$) are not precisely known, and thus it is difficult to compute $\bm{C}$. The concept of \textit{structural controllability} was proposed by C. T. Lin to overcome this limitation \citep{L1974}, which is a property for linear systems slightly weaker than \textit{complete controllability}, as given by Kalman's condition. Yet in practice, whereas complete controllability adds to mathematical elegance but barely provides practical utilities, structural controllability is almost always sufficient to guarantee feasible control of dynamic systems. The theory adopts the abstract graph representation of systems, with nodes denoting state variables, and directed links denoting interdependencies between variables. Based on the graphic representation, Lin defined the basic concepts and structures in dynamic systems that are structurally controllable/uncontrollable (Figure 2): 

\textbf{Theorem 0} \textit{Structural Controllability} 

Under the graphic view of linear dynamic systems, a node is \textit{non-accessible}, if and only if it cannot be reached by any control input (Figure 2d). In the directed system graph, if there exists a set of nodes, denoted by $S$, which contains $k$ nodes ($|S| = k$), such that there are no more than $k-1$ nodes in the companion set $T(S)$, which is the set of all the nodes that have direct edges connected to any of the $k$ nodes in the original node set $S$, then we say the graph contains a \textit{dilation} (Figure 2d). Systems which contain non-accessible nodes or dilations are \textit{uncontrollable}. In the system graph, a \textit{stem} (Figure 2a) or a \textit{bud} (Figure 2b) does not contain any non-accessible node or dilation, so they are two end-member structures that are structurally controllable, and sophisticated structures built up uniquely from these two structures (\textit{cacti}, Figure 2e) are always structurally controllable. Note that structural controllability requires that all edges in the graph (elements in matrix $A$) indicate free parameters in real space; the system is not guaranteed to be structurally controllable if certain interdependencies between state variables only have a limited parameter range.$\square$

\begin{figure}[!ht]
\centering   
	\includegraphics[width=6in]{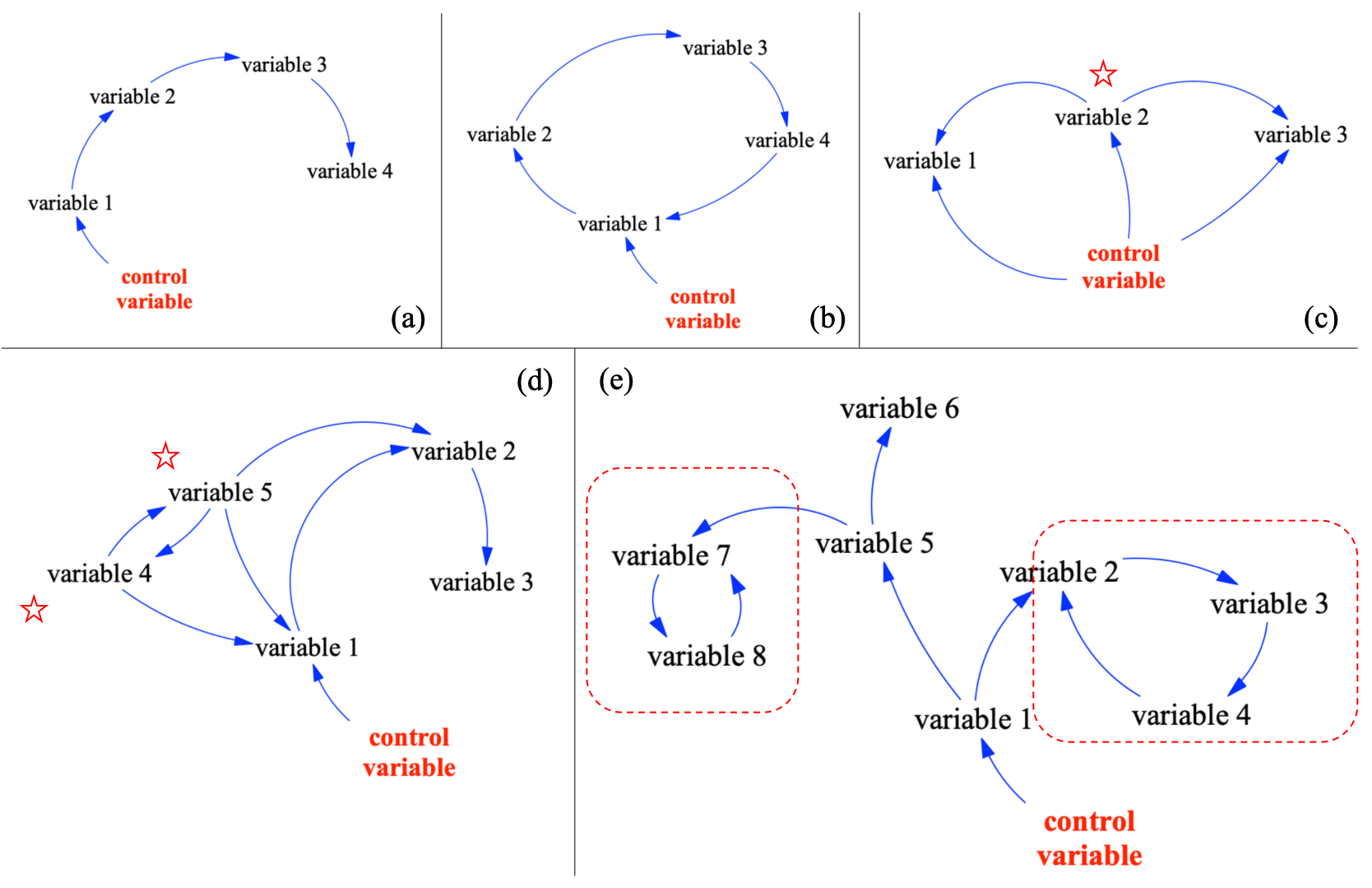}  
\caption{End-member structures of linear dynamic systems that have clear structural control properties \citep{L1974}. (a) stem; (b) bud; (c) dilation (on variable 2); (d) non-accessible nodes (variable 4 and 5); (e) cactus (a stem connecting two buds in dashed boxes). (a), (b) and (e) are structurally controllable structures; (c) and (d) are uncontrollable.} 

\end{figure}

Note that controllability does not mean that a reached state can be \textit{maintained}, merely that any state can be \textit{reached}. It is for this reason that one could control an entire \textit{stem} or \textit{bud} with a single control input. Apparently, more control efforts are needed to let the system maintain a reached state.

One important corollary of the above theorem is the role of self-loops in helping enhance system control. By adding new self-loops to dilated nodes in an uncontrollable system (e.g. onto variable 2 in Figure 2c), the companion set $T(S)$ could be expanded, since the dilated nodes with new self-loops now belong to this set, and as long as the companion set $T$ is as large as the original set $S$, the dilation will disappear and the system will become controllable. It is a widely known rule of thumb in modeling practice that self-loops on stock variables (first-order exit control) are desired in SD models; the above corollary may serve as a theoretical justification for this heuristic.

Structural control principles are powerful tools for determining the controllability of complex systems, since they could also work on time-variant non-deterministic systems, where the dependency matrix $\bm{A}$ is unfixed and only partially known: zero entries are fixed, nonzero entries could take different values, i.e., the interdependency relationships between state variables always hold but values may evolve. The graphic criterion is also more straightforward and far less computationally costly than Kalman's rank condition. 

Relying on the concept of structural controllability, recent studies have extensively addressed the structural control principles of complex networks and great theoretical progress has been made \citep{Let2011, Cet2012, LB2016}, triggering a number of successful applications on various fields \citep[e.g.][]{NA2013, Ket2016, Vet2016}. However, so far most results apply to linear systems; dynamics in nonlinear systems are far more complicated as the interdependencies between variables take non-additive functional forms. In order to extend structural control theories to nonlinear systems, some studies turn to the linearization of nonlinear dynamics \citep{Wet2015}, while other studies look directly into the specific nonlinear feature of the system \citep{HK1977}. With decent efforts, new concepts on the control principles of nonlinear systems have been developed in recent attempts, such as exact controllability (extending structural controllability to a more general framework) \citep{Yet2013}, and control based on feedback vertex sets \citep{Fet2013,Met2013}. Although much progress has been made and many exciting applications have taken place \citep[e.g.][]{Get2015, Zet2016}, the study of control principles for nonlinear dynamic systems is undoubtedly still at the early stage.

\section*{Controllability of SD Models}

It is a natural idea to apply structural control principles for dynamic systems onto SD models. Theoretically speaking, SD models belong to highly nonlinear time-invariant dynamic systems, where the interdependencies between variables are not time-dependent, but often have a limited value range and are not freely tunable. High nonlinearity of SD models derives from the constant use of table functions, min/max formulations, delays etc., on top of the more common usage of exponentials, logarithms, multiplications and other nonlinear analytical functionals. It is then almost impossible to linearize SD models and therefore the structural control principles for linear systems may not directly apply to SD and need certain translation. Complimenting classic theories, specific control principles for SD models are to be derived, considering unique SD formulations and specific fields of application.

\subsection*{\textit{SD models and physical dynamic systems: critical differences in graphics}}

It is easy to mistake a typical SD model (Figure 1) with a physical dynamic system (Figure 2). Although the abstract graphic views look similar on the two sides, they are different in the underlying dependency relationships. In the diagram of a physical dynamic system, an arrow connecting (state) variable $x_1$ to variable $x_2$ indicates that the velocity in $x_2$ depends on $x_1$, i.e. $\dot{x_2} = c\times x_1$; in a typical SD diagram, an arrow from $x_1$ to $x_2$ indicates a direct dependence, i.e. $x_2 = c\times x_1$. 

More formally, in discrete time notations, a physical dynamic system with state variable $x$ and control input $u$ could be represented by (Figure 3a):

\begin{equation}
\begin{aligned}
\frac{x(t+\delta t) - x(t)}{\delta t} = f'(x(t),x_1(t),...x_n(t)) + b'u(t+\delta t) \\ 
\Rightarrow x(t+\delta t) = f(x(t),x_1(t),...x_n(t)) + bu(t+\delta t) \ \ t>0
\end{aligned}
\end{equation}
where $\delta t$ is the computational interval in discrete time and $f$($f'$ as well) is a linear or nonlinear dependency. On the other hand, a SD model with a similar graphic structure should be represented by (Figure 3b):
\begin{equation}
x(t+\delta t) = f(x_1(t),...x_n(t)) + bu(t+\delta t),\ \ \ \ t>0.
\end{equation}

Comparing (3) and (4), since $x(t+\delta t)$ depends on $x(t)$ in (3) but not in (4), there is a fundamental difference in the graphic representation of SD models and physical dynamic systems. The counterpart of a node (state variable) in physical systems (Figure 2) is a stock variable in SD models (Figure 1), and an arrow in Figure 2 corresponds to a flow in Figure 1. Variables besides stocks in Figure 1 (known as auxiliary variables in SD) are in fact intermediate components of interdependencies between stock/state variables. It is expected that the dependency structures of SD models would be quite simplified when abstracted into the style of Figure 2, since an SD model normally contains few stocks; however, given the huge complexity in systems' nonlinearity, SD models need to be represented in a finer scale in order to address their complex dynamics. Hence in most studies they are, and should be, presented as in Figure 1, more detailed than the graphic representation of physical systems. 

To facilitate our discussion, in this study, we abstract SD models into a graphic view similar to that of physical systems, as adopted by standard structural control analysis \citep[e.g.][]{Let2011}, but with the following adapted conventions: stock variables are represented as filled (blue) nodes, auxiliary variables are represented as empty (blue) nodes; interdependencies between stock variables are represented by arrows (solid or dash; see following sections); control inputs are represented by (red) squares (Figure 3 and on). 

\begin{figure}[!ht]
\centering   
	\includegraphics[width=6in]{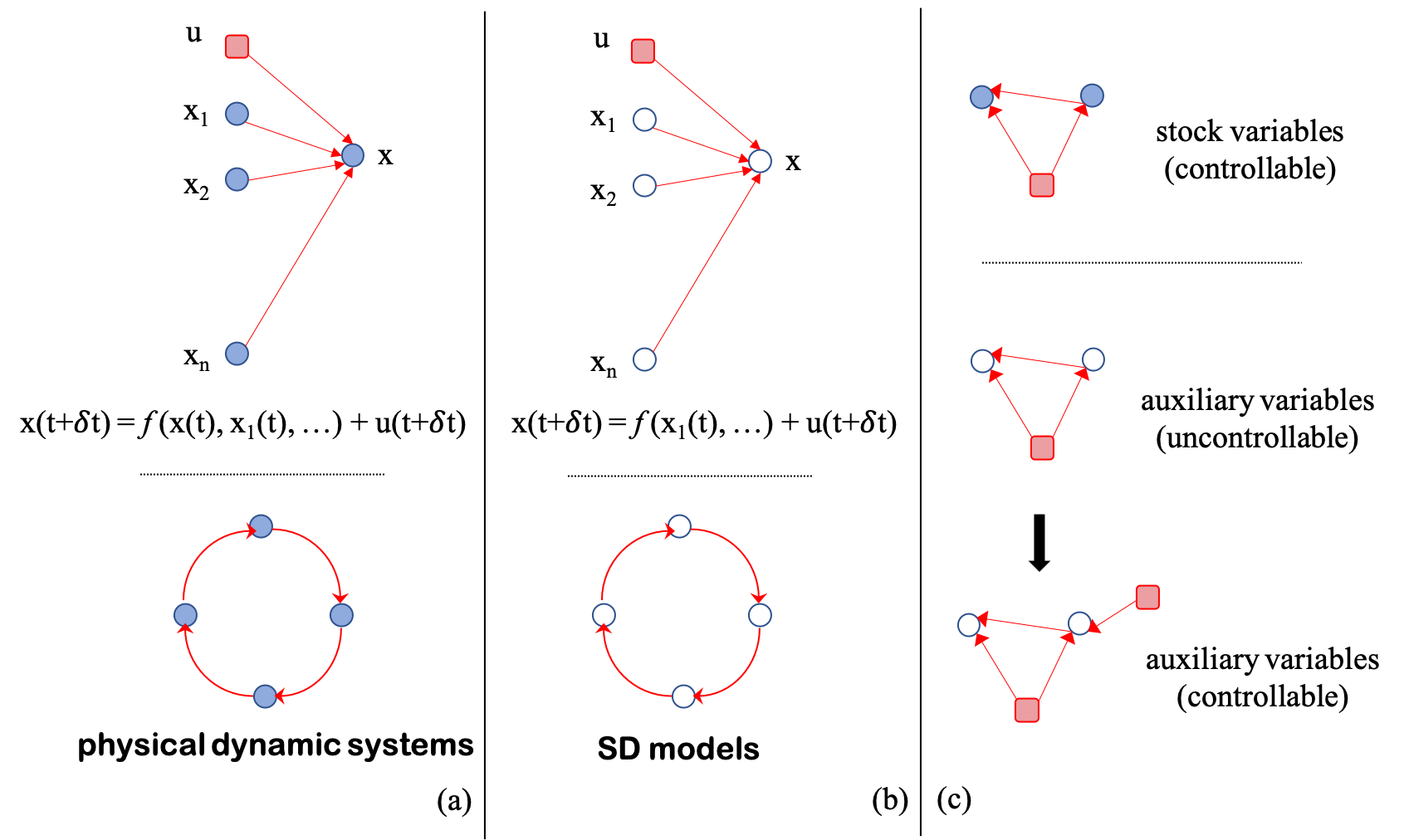}  
\caption{Underlying differences in the similar graphic representations between (a) physical dynamic systems and (b) SD models. (c) Sample structure demonstrating the different control properties of physical dynamic systems and SD models under similar graphic views.} 

\end{figure}

Consider the triangle structure shown in Figure 3c. According to classic structural controllability theory (Theorem 0), if the nodes are stock variables (filled), as are the state variables in physical dynamic systems, the system is controllable since no dilation exists. However, if the nodes are auxiliary variables in SD (empty) and thus the links do not represent the rate of change, then a dilation could be identified and the system is no longer controllable; only when an extra control input is added to the system will it become controllable again. This simple example demonstrates the root difference in the graphic representation between SD models and physical dynamic systems, and echoes our argument that the structural controllability theory does not directly apply to SD models and needs to be adapted for specific SD structures. 

\subsection*{\textit{An important theorem}}

Keeping in mind the above discussion of the difference in the graphic view between physical dynamic systems and SD models, the structural controllability theory (Theorem 0) could be adapted for SD. The two end-structures, the \textit{bud} and the \textit{stem}, frequently appear in SD diagrams, but in most cases incorporating both stock variables (filled nodes) and auxiliary variables (empty nodes). Through further inspection, it is critical to point out the following fact: in an SD structure that only consists of \textit{buds} and \textit{stems} (i.e., \textit{cacti}) , when {replacing some stock variables with auxiliary variables}, the transformed structure on stock variables (omitting all auxiliary variables to follow the conventions in Figure 2) is always \textbf{still} a structure uniquely consisting of \textit{buds} and \textit{stems}, with the possibility that a \textit{bud} of stock variables may degenerate into a \textit{stem}. 

Reflecting on the Theorem 0, this observation yields a significant result:

\textbf{Theorem 1} \textit{Sufficient Condition on Structural Controllability for SD Models} 

Under the graphic representation of SD diagrams consisting of both auxiliary variables (empty nodes) and stock variables (filled nodes), if the structure transformed from the original graph, when replacing all auxiliary variables with stock variables, is structurally controllable (meaning that it is a dynamic system consisting of only \textit{buds} and \textit{stems}, concerning stock variables), then the original SD model is structurally controllable. $\square$

This theorem expresses a sufficient condition for determining the structural controllability of SD models, utilizing the fact that SD diagrams contain two distinct categories of variables. It is a powerful tool for the structural control analysis (SCA) of SD models (see following sections).

\subsection*{\textit{Identification of control inputs}}


In current SD terminology, variables are categorized as stock (level) variables and auxiliary variables. Stock variables have long been paid great attention to by SD modelers and are easy to pin down, yet it is often neglected by many that, within the category of auxiliary variables, there are different subgroups of variables that play distinct roles in the system. Specifically, one need to distinguish \textit{parameters} from \textit{control inputs}. Both parameters and control inputs have no dependence within the model boundary and use exogenous values; however, they play different roles in determining system's dynamics: at steady state, manipulating control inputs changes the levels of stock variables, whereas tuning parameters does not influence the steady-state values of stocks. Similar to the situation in physical systems, in order to conduct control analysis on SD models, the identification of control inputs is a necessary first step. 

For example, in the stock management structure (Figure 7), which is a classic SD model component, both \textit{Desired Supply Line} (SL*) and \textit{Supply Line Adjustment Time} (SLAT) are auxiliary variables that require exogenous inputs (a similar pair for \textit{Desired Stock} and \textit{Stock Adjustment Time}). In current SD terminologies, the two variables are not distinguished by straightforward notations. However, experienced SD modelers would recognize that, \textit{Desired Supply Line} is a control input, while \textit{Supply Line Adjustment Time} is a parameter. Empirical knowledge may be called on here, as \textit{Desired} values influence stocks and any \textit{Adjustment Time} is a granted parameter. Similarly, variables named with \textit{sensitivity}, \textit{threshold}, \textit{lag}, \textit{delay time}, etc. are also known to be parameters in almost all cases. Nevertheless, although modeling experiences would be sufficient for the separation of control inputs and parameters under many circumstances, a more rigorous discussion on this task is certainly welcomed. With such an idea in mind, we make the following arguments, which immediately derive from basic control theory (equation 1) yet should be brought to special attention. 

Assume an auxiliary variable $z$ from an SD model, which has no dependence on system's internal components and requires exogenous input. The following property determines whether $z$ is a parameter (P) or a control input (C):
\begin{equation}
PorC(z) \triangleq \sum_{i,j} PorC_{i,j}(z) = \sum_{i,j} \frac{\partial}{\partial x_j} (\frac{\partial \dot{x_i}}{\partial z}) = \sum_{i,j} \frac{\partial^2 \dot{x_i}}{\partial z \partial x_j}\ \ for\ \forall i,j = 1,2,...,N,
\end{equation}

where $\bm{x}$ represent stock variables in the system. Taking derivatives on equation (1), one could see that, if $z$ is a control input, i.e., an element in $\bm{u}$, then $PorC_{i,j}(z)=0$ for any $i,j$, so $PorC(z)=0$. On the other hand, if $z$ is a parameter of the model, i.e. an entry in matrix $\bm{A}$, then $PorC_{i,j}(z)\neq 0$ for some $i,j$, so $PorC(z)\neq 0$. 

This $PorC$ property then rigorously distinguishes parameters from control inputs, among exogenous auxiliary variables in SD models. It is by no means a novel invention, yet this simple property is useful for the specific task of control inputs identification, which is the first step in the structural control analysis on SD models (see following sessions).

\subsection*{\textit{Levels of system control in SD models}}

After the successful identification of control inputs, we are now able to establish specific control principles for SD models. In SD simulations, we expect endogenous model variables to follow desired trajectories; based on our expectations of how model variables could be manipulated by control inputs, different levels of system control are to be realized in SD models. 

At the lowest level of manipulation, the value of an endogenous variable could be at least slightly changed by control inputs, i.e., its dynamics are not completely independent of system control. This level of control demonstrates \textit{accessibility}.

\textbf{Definition}  \textbf{\textit{Accessibility}} A variable in an SD model is \textit{accessible}, if and only if its value could be (at least) changed by certain control inputs. It is \textit{non-accessible} if its dynamics will not be influenced at all by any control input. An SD model is \textit{accessible} if it does not contain any \textit{non-accessible} variable. $\square$

In theory, both stock variables and auxiliary variables could be non-accessible (Figure 4a), although most SD models expect stock variables to be accessible. As mentioned, unlike physical systems, non-accessibility is often \textit{not} an undesired property for dynamic models in social sciences. In many situations, a simulation model for real-world problems is built to demonstrate certain dynamics that appear and develop in a natural way, in which case control inputs are excluded. In other cases, specific parts of a model are expected to be non-accessible and follow natural dynamics. It is not always the case that SD modelers expect the whole system to be entirely \textit{accessible}, where control inputs could influence all endogenous variables. 

A higher level of control provides modelers with enhanced manipulation over the system. In physical dynamic systems, the realization of system control requires that the interdependencies between variables (elements in $A$) could be freely manipulated in real space $\mathbb{R}$; the system may deviate from being structurally controllable if any interdependency only has a limited range for manipulation. However, in SD models, where interdependencies are represented by intermediate auxiliary variables, the free-tuning of variables in the real space is often not realized. The prevailing use of nonlinear formulations in SD models such as MIN/MAX cutoffs, table functions, and various analytical functionals greatly dampens this task and results in narrowed value ranges. Taking note of this discrepancy between physical systems and SD models, we define an advanced level of system control especially for SD variables, on top of \textit{accessibility}: we call an (accessible) variable \textit{spanning}, if it could span the entire value range in simulations by control inputs, instead of residing in a narrowed value range, in which case it is \textit{non-spanning} (Figure 4b).

\textbf{Definition} \textbf{\textit{Spanningness/Non-spanningness}} A variable in SD models is \textit{spanning} if and only if its entire value range (in theory, the real space) could be reached with the help of control inputs. It is \textit{non-spanning} if certain values within its value range cannot be reached. More importantly, a mathematical formulation for the interdependency between variables is a \textit{non-spanning} formulation if it results in a \textit{non-spanning} variable. $\square$

Similar to non-accessibility, non-spanningness is \textit{not} an avoided property in simulation models. In fact, many modelers strongly desire that certain variables in their models do not span the entire theoretical value range but instead only fluctuate within a limited range, which may then represent specific insights of their model. Analytical nonlinear formulations and manual cutoffs (MIN/MAX, table functions) are used in SD models exactly for this purpose. To a great extent, non-spanning formulations largely help retain operational interpretability of model analysis by downgrading the controllability of the system.

Finally, some models may achieve the highest level of controllability, as the established concept in classic control theory. Note that in SD practice people mainly focus on stock variables, hence the controllability of stock variables and the controllability of the whole system are to be differentiated. However, the latter is rarely achieved since it is only theoretically possible to control all auxiliary variables in an SD model.

\textbf{Definition} \textbf{\textit{Controllability/Quasi-controllability}} An SD model is controllable/quasi-controllable, if and only if any start state of the system (state containing the entire panel of variables/state containing only stock variables) $\bm{X}(t_0)$ could be driven to any possible end state $\bm{X}(t_1)$ in finite time through control inputs $\bm{U}$, i.e., mathematically: 
\begin{equation}
\bm{X}(t_0) \xrightarrow{\bm{U}}\bm{X}(t_1)\ \ \ \ \ \forall \bm{X}(t_0),\bm{X}(t_1) \in \bm{X}, t_1 - t_0 < + \infty
\end{equation}
$\square$ 

\begin{figure}[!ht]
\centering   
	\includegraphics[width=4in]{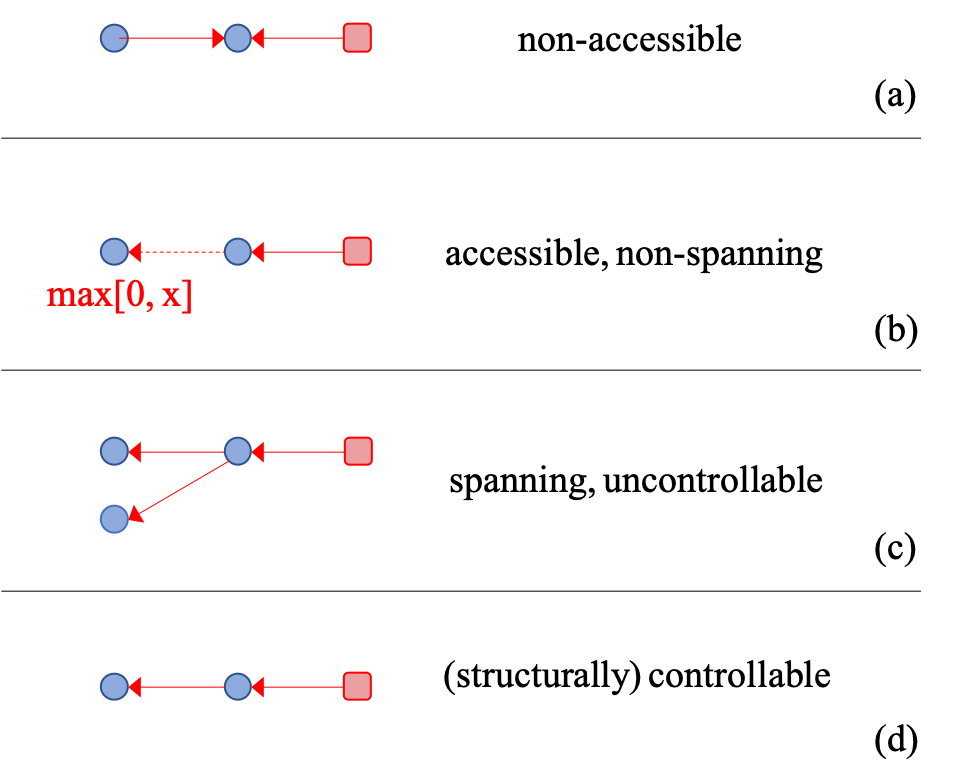}  
\caption{Multiple levels of system control in SD models. (a) System with non-accessible nodes (the leftmost node). (b) Accessible but non-spanning nodes (the leftmost node). (c) Spanning but uncontrollable system (with a dilation). (d) Structurally controllable system.} 

\end{figure}

So far we have conceptualized the three-level system control in SD models (Figure 4). For stock variables, our definitions of \textit{accessibility} and \textit{controllability} recover the classical definitions in physical dynamic systems. The second level of system control (\textit{spanningness}) does not appear in physical systems and is the special property for the auxiliary variables in SD models. We consider \textit{spanningness/non-spanningness} as an important concept in SD, since in SD models a great portion of variables are auxiliary variables and extensive nonlinear relationships are used (see the following section). For simplicity, we no longer distinguish \textit{controllability} and \textit{quasi-controllability} in this study, as our analysis mainly focuses on the controllability of stock variables in SD.

\section*{Control Properties of SD Model Components}

So far we have established the structural control principles for SD models. Taking note of the uniqueness of SD model structures among dynamic systems, the general structural controllability theory (Theorem 0) has an extended discussion in SD (Theorem 1). The distinction between auxiliary variables and stock variables underlies the three-level system control in SD structures. Under this theoretical framework, we are able to study the control properties of important SD model components, including structures on stock variables (single stocks, co-flows, chains of stocks, etc.), manual nonlinear formulations (MIN/MAX and table functions), feedback loops, and time delays. For decades, heuristics in model building have been extensively used in SD practice; it is shown that, to a certain extent, some of these modeling heuristics may be readily justified by the structural control analysis of SD models. 

\subsection*{\textit{Structures on stock variables}}

As mentioned, stock variables are a major focus in most SD models. Three structures on stock variables are often seen in SD models: single stocks, chain of stocks, and co-flows; following the graphic conventions of this study, their structures are transformed into abstract views (Figure 5). In typical SD diagrams, the net rate of change in stocks is separated into enter rate (inflow) and exit rate (outflow); however, the dependencies of stock variables essentially represent their net rate of change. This means that a dependency arrow connected to a node in the graphic view corresponds to either the enter rate or the exit rate of the equivalent stock variable in the SD diagram. Without loss of generality, we uniformly match the dependencies arrows in abstract graphs to the exit rates of stocks in SD diagrams (Figure 5a-5c). 

The control properties of these structures on stocks could be readily identified. Since they uniquely consist of stocks, the general structural controllability theory (Theorem 0) directly applies. Specifically, we have the following propositions:

\textbf{Proposition 1} In SD modeling practice, first-order exit control on stock variables enhances the controllability of the models. $\square$

Proposition 1 holds true because first-order exit control on stocks essentially corresponds to self-loops on the nodes. As mentioned before, self-loops help eliminate dilations in the system and thus help enhance system control. First-order exit control on stocks has long been considered as good modeling practice in SD, since building stocks without first order control will likely make the stocks become negative. This heuristic might now be justified by structural control arguments. 

\textbf{Proposition 2} Co-flows and chains of stocks are structurally controllable, as long as the first stock in the series is controllable. $\square$

According to Theorem 0, co-flows and chains of flows (e.g., aging chains) are \textit{stems}, which ensure structural controllability. Moreover, adding more stocks in co-flows and aging chains does not change the controllability of the system. As is known to experienced SD modelers, it is often impossible that a downstream stock variable in a chain of stocks or co-flows could fall out of control while the upstream stock is being successfully manipulated.

\begin{figure}[!ht]
\centering   
	\includegraphics[width=4in]{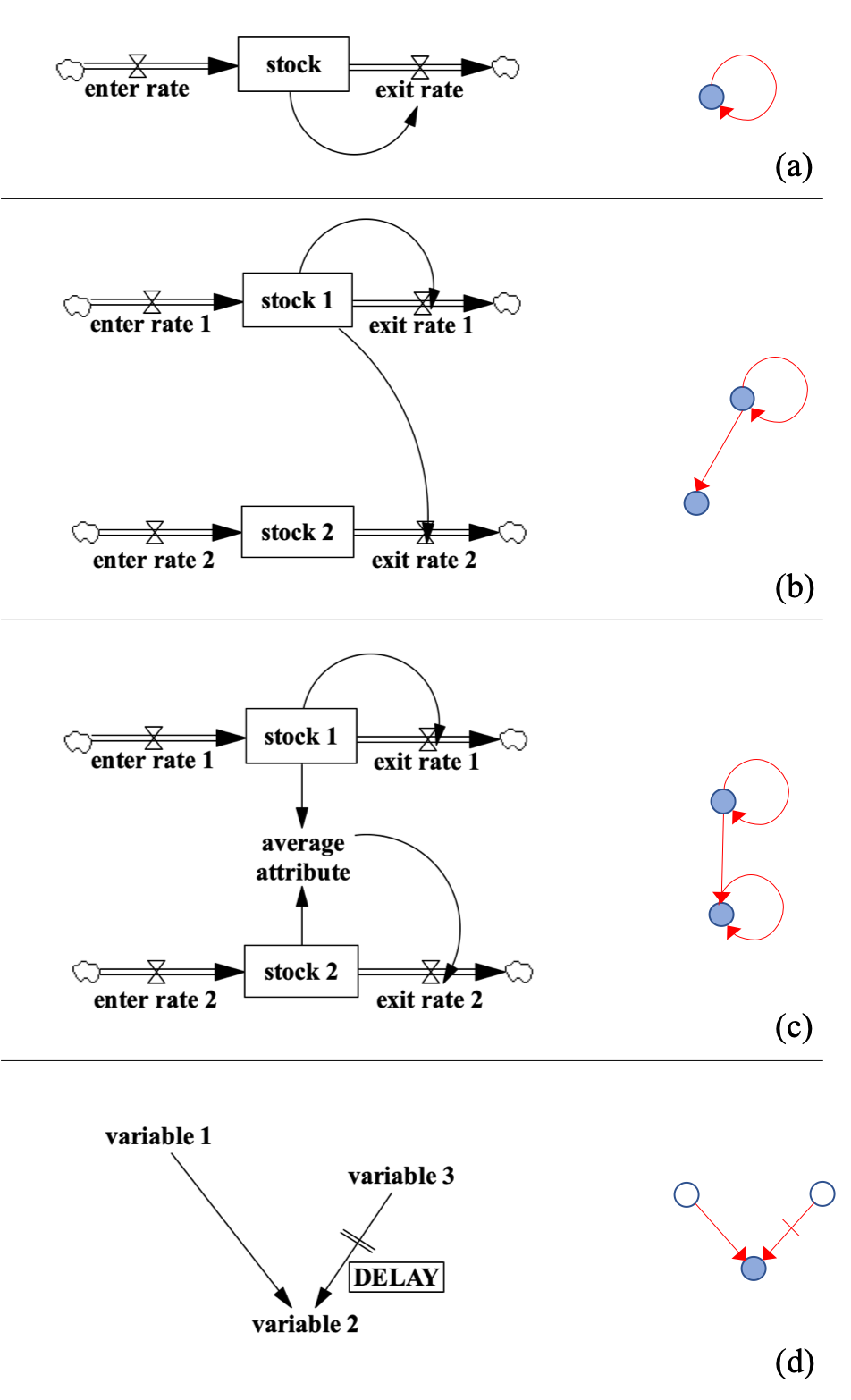}  
\caption{Control properties of important SD model components. The SD diagram and the abstract graphic view are shown for each structure. (a) Single stock.  (b) Chain of stocks. (c) Co-flow. (d) Delay on non-stock (auxiliary) variables.} 

\end{figure}

\subsection*{\textit{Non-spanningness and nonlinear SD components}}

The prevailing nonlinearity in variables' interdependencies determines a crucial distinction between SD models and physical dynamic systems. Unlike in power grids or transportation networks where many sensors, controllers, stations, etc. constitute the system, in social, economic and managerial systems, variables are often of different nature; thus many complex interdependent relationships may be built in these models. Nonlinear functionals such as exponentials, powers, logarithms are often used, on top of simple multiplication, although the usage of more complicated functionals that are hard to interpret is also largely avoided in SD practice. 

In cases where the relative rigid curvatures of analytical functionals are not desirable, SD modelers turn to manual nonlinear relationships such as MIN/MAX and table functions. These formulations transform inputs to desired outputs and inevitably narrow the value range of the output variable, in most cases an auxiliary variable. According to our established control principles, this naturally introduces \textit{non-spanning} relationships between input and output variables. The value of output variables may be constrained in the half space (for MIN/MAX, exponentials, powers of 2, etc.) or even in a closed interval (for table functions). To represent such non-spanning relationships, we make an important point in our graphic conventions: we denote links whose end nodes embody non-spanning relationships of start nodes by dash lines, and links with no inclusion of non-spanning formulations as solid lines.

As has been discussed, non-spanningness is often a desired property in SD models. The realization of structural controllability relies on the free-tuning of interdependencies between stock variables; any non-spanning formulation dampens the free-tuning of variables and thus limits the system's controllability. Since achieving complete controllability of the system is often not what simulation models are built for, it is appropriate to argue that non-spanning formulations are indispensable for SD models, and perhaps for other social and economic dynamic models as well. Specifically, the role of table functions is to be highlighted, since they often constrain the value range of variables that is available for system control to the largest extent. 

\subsection*{\textit{Feedback loops: reinforcing and balancing}}

Next we discuss the control properties on a major SD component, the feedback loops. Unlike other SD model components, feedback loop is a more general and inclusive concept. The structure of a feedback loop is complex, probably consists of both auxiliary variables and stock variables, and both linear and nonlinear/non-spanning relationships. As a result, one should note that feedback loops in SD settings are different from typical feedback structures in physical dynamic systems, where all the nodes are exclusively state variables, and in most cases no non-spanning relationship appears. Therefore, the control theory for dynamic systems based on Feedback Vertex Set (FVS) \citep{Fet2013, Met2013} in general does not apply to SD models.

Given the vast complexity of feedback loops in SD models, results on their control properties are difficult to pin down. It may nevertheless be useful to distinguish the spanningness of reinforcing feedback loops and balancing feedback loops. One observation is that, reinforcing feedback loops should be spanning, if they are going to grow without limit. If one variable within the loop is not spanning and the value is capped, the unstoppable growth will not be triggered. On the opposite, as one primary utility of balancing feedbacks is to stabilize the system, variables in balancing feedback loops are expected to fall around certain values, i.e., they are often not spanning. Balancing feedbacks help constrain the dynamics of the model, as is acknowledged by many SD practitioners, and in this sense, help limit the controllability of the system. In order to achieve the spanning of variables in a balancing feedback loop, one almost always needs delays (except for the situation where the variable is single-valued, i.e., it is always spanning). Delays desynchronize the dynamics of different variables in the loop and essentially make their trajectories fluctuate. It is well known that delays are necessary in balancing feedbacks in order to trigger oscillatory dynamics; once again, one might justify such a heuristic in modeling practice from the structural control perspective. To summarize, we have the following proposition. Note that we enlarge the concept of spanningness by using it to describe feedback loops.

\textbf{Proposition 3}	 Reinforcing feedback loops are spanning, if all variables in the loop could grow without limit. Balancing feedback loops are almost always not spanning if no delay exists in the loop. $\square$

\subsection*{\textit{The special role of delays}}

It is widely acknowledged in the field that time delays (and smooths) are important components in SD models. The existence of any delay implicitly associates with a certain stock; the input variable goes into the imaginary stock, and the output rate of the stock is modulated by delay time and order. For structural control analysis, it is then important to explicitly represent the stock. Hence, in our graphic conventions, the variables at the end of delay arrows (arrows with a vertical bar in the middle) are viewed as stock variables (Figure 5d). 

Failing to identify and represent delays in the model may result in false structural control analysis. Consider the graphic structure in Figure 6. If the delay is not highlighted and the variable at the end of the delay arrow is viewed as an auxiliary variable, the system would be determined as structurally controllable (Figure 6a); however, once the underlying stock beneath the delay is correctly represented, a dilation appears and the system is not controllable any more (Figure 6b). Structures like this example are not hard to find in SD models, and time delays should always be explicitly represented.

\begin{figure}[!ht]
\centering   
	\includegraphics[width=4in]{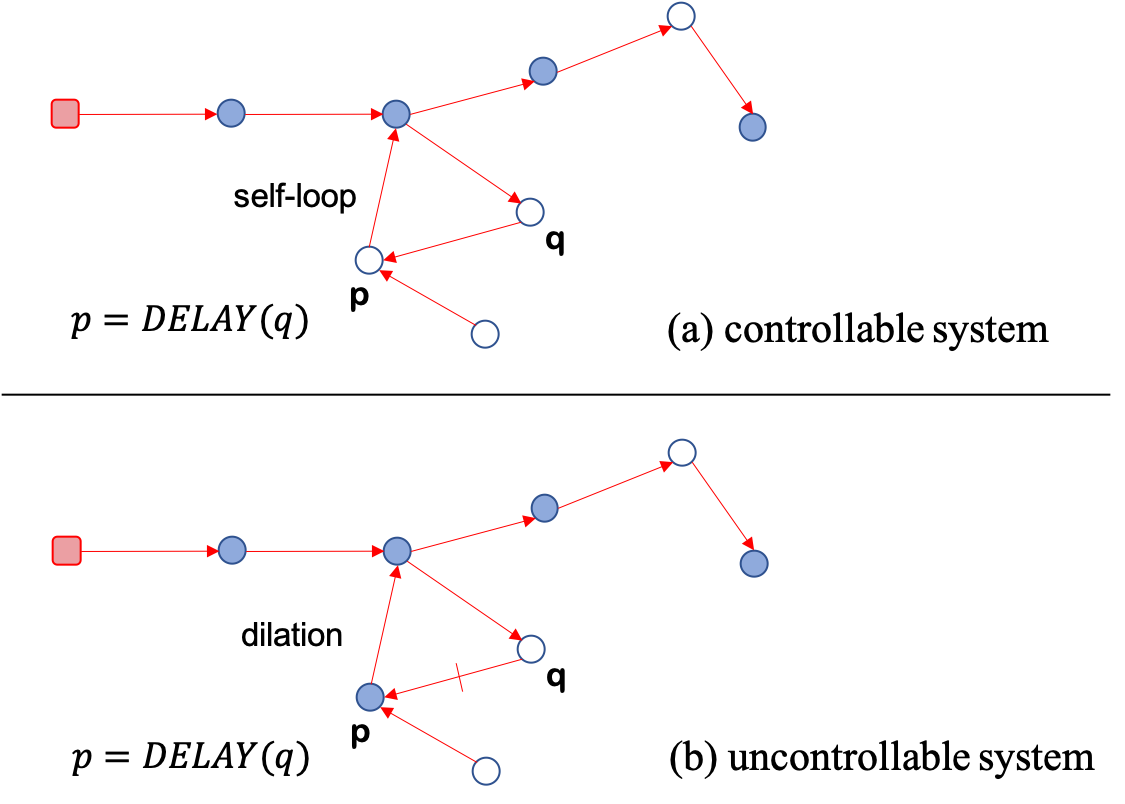}  
\caption{Special concerns on time delays in SD models. Delays are associated with implicit stock variables; failing to represent these hidden stocks may lead to incorrect structural control analysis of the model. A delay from variable $q$ to $p$ exists in the sample structure. (a) Incorrect representation: the delay is not identified; system is structurally controllable with a self-loop. (b) Correct representation: delay is identified and the implicit stock is captured. The system is actually uncontrollable; the self-loop is in fact a dilation.} 

\end{figure}

One subtlety about delays is that, unlike explicit stock variables which are assumed to always take positive values, the stocks underlying time delays could possibly take negative values. Hence, the bi-direction values of the underlying stocks intrinsically let delays be able to serve as stabilizers of the system. Indeed, normal stock variables are not stabilizing because they do not allow exit rate to exceed enter rate, either for a long time or to a large extent; not emitting enough materials from stocks fails to bring the system back to order and, instead, accumulates hazards for a potential system collapse. This discussion responds to the utility of time delays and smooth functions, which are extensively used in SD models to slow down the dynamics and filter out high frequency fluctuations. Once again, control principles lie in the core of such modeling heuristics.

\section*{Example: Structural Control Analysis On the Stock Management Structure}

Under the theoretical framework of multi-level structural control in SD models that we have established in this study, and relying on the discussion of the control properties of typical SD model components, we are able to conduct structural control analysis (SCA) on SD models. As an example, in this section we show a sample analysis on the generic Stock Management structure \citep{S2000}, which is a classic SD model structure in supply chain management (Figure 7a). The theoretical and graphic tools that we developed are brought together to study the model.

The structure consists of 2 stock variables (\textit{Supply Line}, \textit{Stock}) and 6 auxiliary variables (\textit{Order Rate}, \textit{Indicated Orders}, \textit{Adjustment for Supply Line}, \textit{Adjustment for Stock}, \textit{Desired Acquisition Rate}, \textit{Expected Loss Rate}). Four exogenous variables remain in the model diagram (\textit{Desired Supply Line (SL*)}, \textit{Supply Line Adjustment Time (SLAT)}, \textit{Desired Stock (S*)} and \textit{Stock Adjustment Time (SAT)}) and are candidates for control inputs. 

According to equation (5), the $PorC$ property is calculated for each of the four candidate variables. The dynamics of \textit{Supply Line (SL)} is governed by $\dot{SL} \propto (SL^*-SL)/SLAT$, so one could determine that $PorC(SL^*) = 0$, and $PorC(SLAT) = SLAT^{-2}\neq 0$. Hence \textit{Desired Supply Line (SL*)} is a control input, and \textit{Supply Line Adjustment Time (SLAT)} is not a control input but a parameter. Similarly, \textit{Desired Stock (S*)} is another control input, and \textit{Stock Adjustment Time (SAT)} is a parameter. 

\begin{figure}[!ht]
\centering   
	\includegraphics[width=6in]{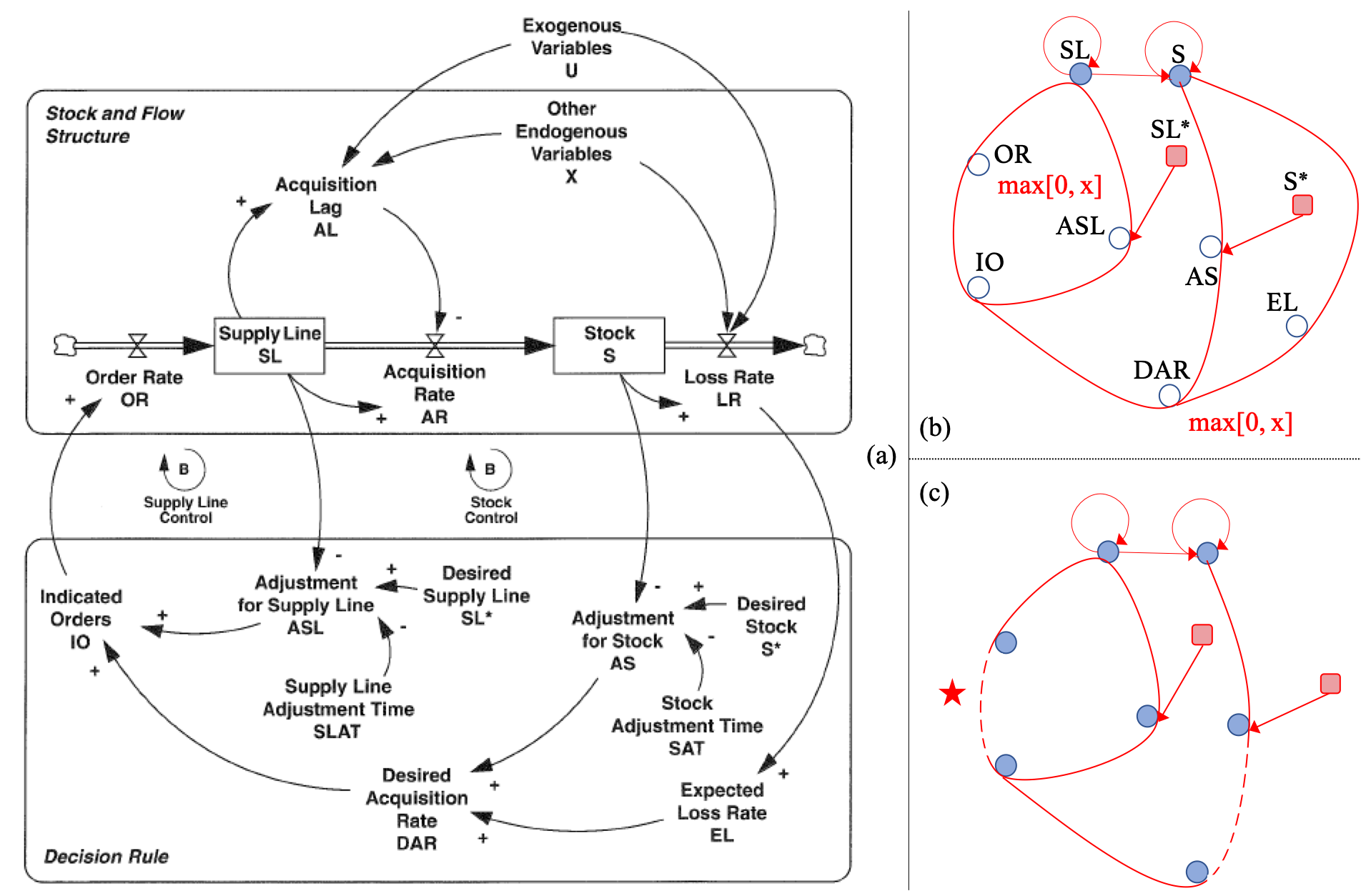}  
\caption{Structural control analysis on the stock management structure. (a) Original SD model (adapted from Figure 17-6 in \citet{S2000}). (b) Abstract graphic view at the first stage. (c) Modified model representation for structural control analysis. Note $U$ and $X$ in the original SD model are not system variables and are omitted in (b) and (c). The two pathways going from $S$ to $DAR$ in (b) are aggregated in (c).} 

\end{figure}

Incorporating the two identified control inputs in the graph and omitting the parameters, we formulate the graphic representation of the stock management structure, consisting of 10 nodes (Figure 7b). The mathematical formulations of the model are inspected. There are two MAX functions appearing in the system, in the formulation of \textit{DAR} and \textit{OR}, respectively, both of which are used to ensure that the enter rate into the (chain of) two stocks $SL$ and $S$ is non-negative, i.e., there is no negative order of products. These two MAX functions indicate two non-spanning relationships, which are denoted by dashed lines. There is no explicit delay in the structure, so no auxiliary variable (empty node) is switched into the stock variable representation (filled node). 

After pinning down the graphic view of the model structure, now we are able to study the control properties of the system. Both of the two stock variables have first-order exit control to enhance system controllability (Proposition 1). The sufficient condition (Theorem 1) for the original structural control theory (Theorem 0) is adopted. We replace auxiliary variables with stock variables while maintaining the dependency network, with the two pathways going from $S$ to $DAR$ aggregated. The resulting structure is essentially the overlap of two buds (Figure 7c), which could be structurally controlled by any of the two control inputs ($SL^*, S^*$), if the two non-spanning relationships are dismissed, i.e., dash lines are recovered by solid lines. Supposing this is the case, then since the modified structure is structurally controllable, the original SD model is also structurally controllable over the two stock variables by each of the two control inputs (Theorem 1).

Intuitively, our identification and analysis is consistent with the functioning of the system: both $SL^*$ and $S^*$ are exactly the desired values (as in their names) that one sets for the two stocks, and the stock management structure is exactly built with the expectation that the model dynamics would bring the stocks to their target levels. Nevertheless, structural control analysis pins down additional theoretical consolidation for the model, and injects constructive details into the discussion.

However, the two non-spanning relationships are not dispensable. The two MAX functions ensure that stock $S$ and supply line $SL$ do not become negative, which is a critical operational constraint for stock management. Indeed, the non-spanning relationships deliberately prevent the system from being completely structurally controllable, and only let the model remain partially structurally controllable to ensure positive stocks. In this case, as expected, the operational interpretability of the model is successfully enhanced by compromising the system's controllability to a certain extent.

It would be useful to go a further step in the analysis and imagine an alternative scenario. Suppose that the MAX function on $OR$ is dismissed and the one on $DAR$ is retained, i.e. the dash line near the star mark becomes solid. This might correspond to a situation where the items in the supply line could be retrieved. Now, the roles of the two control inputs, which are roughly the same in the original analysis, are distinguished in the new situation. While \textit{Desired Stock (S*)} still demonstrates partial controllability due to the remaining dash line, \textit{Desired Supply Line (SL*)} is now a control input that has full structural controllability over the system. This is consistent with our assumption, since on a more flexible (retrievable) supply line, the decision from the supply line side (\textit{Desired Supply Line}) could essentially control the stock side, but not vice versa. It implies that, in this imaginary situation, \textit{Desired Supply Line} is a more important control input than \textit{Desired Stock}, which then deserves more attention in, for example, model calibration or policy analysis. Although not to a full-fledged extent, this discussion nevertheless demonstrates the significant potential of using structural control analysis in studying the relative importance of control variables in SD models, and in providing further guidance in model calibration and applications.

\subsection*{\textit{Procedures for conducting Structural Control Analysis (SCA) on SD models}}

Going over the above sample analysis, the procedures for conducting structural control analysis (SCA) on SD models could be summarized as follows. 

Given an SD model:

\hspace{20pt} \textbf{Step 1} \textit{Identify control inputs of the model using the $PorC$ property.} 

\hspace{65pt} Identify model parameters.

\hspace{20pt} \textbf{Step 2} \textit{Represent the model in the abstract graphic view using established conventions.} 

\hspace{65pt} Distinguish stock variables and auxiliary variables in the representation. 
	
\hspace{20pt} \textbf{Step 3} \textit{Inspect model formulations and modify the graphic representation.} 

\hspace{65pt} Identify and represent non-spanning formulations. 

\hspace{65pt} Explicitly represent delays by switching the graphic notation on certain auxiliary variables.

\hspace{20pt} \textbf{Step 4} \textit{Examine the structural control property \textbf{for each identified control input}.} 

\hspace{65pt} Apply the classical theorem (Theorem 0) and the sufficient condition (Theorem 1) in analysis.

\hspace{20pt} \textbf{Step 5} \textit{Examine the structural control conditions of the model.} 

\hspace{65pt} Discuss control properties of model components and the aggregated system control conditions.

\hspace{65pt} Discuss the roles of different control inputs and corresponding policy implications.

\section*{Concluding Remarks}

The controllability of dynamic systems is a core concept in modern control theory. Methodologically, SD is the field in social and economic sciences that provides modeling solutions to real-world problems through the dynamic system perspective. SD highlights the interdependencies between model variables, and pays great attention to the structural features of the system that give rise to dynamic model behaviors. Studying the structural control principles of SD models, and broadly speaking, of all simulation models that address real-world issues from the dynamic system point of view, is an important theoretical input to social and economic sciences. 

Adapted from the classic structural control theory with specific theoretical establishments, the multi-level control principles for SD models developed in this study contribute a useful and minimally viable tool to SD model evaluation and analysis: the Structural Control Analysis (SCA). Following the stepwise workflow that utilizes the theoretical set-ups as well as the graphical conventions developed in this study, SCA provides multiple functionalities for SD modelers. Theoretically, it conceptualizes different levels of system control in SD, helps clarify control properties of important SD model components, and helps justify various tested heuristics in SD modeling practice. Practically, SCA facilitates the determination and discussion of different variables' roles in the model, based on which it may further provide well-grounded implications for model calibration and policy/strategy analysis.

As has been repeatedly mentioned in this paper, structural control principles for SD models provide useful insights from a perspective orthogonal to that of physical dynamic systems. Unlike physical networks, models of human systems may demonstrate various functionalities in different settings. Typically, simulation models in social sciences may embody descriptive (illustration of observed dynamics), predictive (forecast of future outcomes), and prescriptive (policy and strategic interventions) powers, and system control plays different roles in those venues. Generally speaking, it would be highly desired in modeling social and economic processes if one can show that the model has \textit{no} or \textit{few} control inputs, and that those control inputs only have \textit{limited} impacts to the model. In other words, a full-fledged controllability of the model is in most cases an \textit{obviated} property in human system modeling. Reconciling classic structural control theories and SD conventions with abundant practical details, our establishments in this study take great care of such rooted dichotomies between physical and human systems. The development of the SCA workflow is targeted to such concerns, which aims to help enhance the interpretability of human system models through the discussion on system controllability.

Given their highly nonlinear nature, the control principles for real-world systems are far more complex than those of physical systems. As a result, to release their power in socio-economic applications, well-established classic control theories for linear physical systems, to which the structural controllability theory belongs, need to be translated and complemented to a great extent to deal with the specific nonlinear structures and functionals in social sciences simulation models. The current discussion on the structural control analysis of SD models makes a first and feasible attempt, but is far from being sufficient; a great amount of theoretical, engineering, and applicational effort is expected in this direction, which in a broader sense goes beyond SD and extends to general dynamic system models of all categories in social-economic sciences.

\section*{Acknowledgements}

The author thanks Prof. Rogelio Oliva at Mays Business School, Texas A$\&$M University, whose detailed comments on early drafts greatly helped improve the quality of the paper, and Prof. Hazhir Rahmandad at MIT Sloan School of Management, whose insights helped clarify the standpoints of the discussion of this study. This research is funded by the System Dynamics Group, MIT Sloan School of Management.

\end{document}